\documentclass[10pt,twocolumn,twoside]{IEEEtran}
\usepackage{graphicx}
\usepackage{epsfig}

\usepackage{amssymb,amsmath,amsthm, amsbsy}

\usepackage[UKenglish]{babel}

\usepackage{algorithm} 
\usepackage{algpseudocode}

\usepackage{mathtools}
\usepackage{array}

\usepackage{changepage} 

\usepackage{enumitem}

\usepackage{graphicx}
\graphicspath{{figuresPartII/}} 
\usepackage[export]{adjustbox}

\usepackage{siunitx}

\usepackage{cases}

\usepackage[usenames, dvipsnames]{color}

\usepackage{xcolor}

\definecolor{line_gray}{gray}{0.75}
\definecolor{cell_gray}{gray}{0.9}
\definecolor{cell_grad_blue1}{HTML}{42A5F5}
\definecolor{cell_grad_blue2}{HTML}{64B5F6}
\definecolor{cell_grad_blue3}{HTML}{90CAF9}
\definecolor{cell_grad_blue4}{HTML}{BBDEFB}

\usepackage{colortbl}

\usepackage{multirow}

\theoremstyle{definition}

\newtheorem{lemma}{Lemma}
\newtheorem{definition}{Definition}
\newtheorem{remark}{Remark}

\newtheorem{proposition}{Proposition}

\usepackage{booktabs}

\newcommand{\lsym}{L_{sym}}

\newcommand{\mati}{\Phi}
\newcommand{\matii}{\Psi}
\newcommand{\espacei}{\mathcal{W}}
\newcommand{\espaceii}{\mathcal{V}}
\newcommand{\evecsi}{W}
\newcommand{\evecsii}{V}
\newcommand{\evali}{\phi}
\newcommand{\evalii}{\psi}

\newcommand{\dimevec}{r}

\newcommand\numberthis{\addtocounter{equation}{1}\tag{\theequation}}

\newcommand{\eqnref}[1]{(\ref{#1})}

\def\assnonzerojtheigengaps{{2}}

\setlist[enumerate]{label*=\arabic*.}

\begin{document}

\title{Extending the Davis--Kahan theorem for comparing eigenvectors  of two symmetric matrices II: Computation and Applications}

\author{J.~F.~Lutzeyer and A.~T.~Walden,~\IEEEmembership{Senior~Member,~IEEE}  
\thanks{
Copyright (c) 2019 IEEE. Personal use of this material is permitted. However, permission to use this material for
any other purposes must be obtained from the IEEE by sending a request to pubs-permissions@ieee.org. 
J.~F.~Lutzeyer and A.~T.~Walden are with the Dept. of Mathematics, Imperial College London, London SW7 2AZ, UK (e-mail: jl7511@imperial.ac.uk and a.walden@imperial.ac.uk)
} }
\IEEEpubid{}
\maketitle

\begin{abstract}
The extended Davis--Kahan theorem makes use of polynomial matrix transformations to produce bounds at least as tight as the standard Davis--Kahan theorem. 
The optimization problem of finding transformation parameters resulting in optimal bounds 
from the extended Davis--Kahan theorem is presented for affine transformations. It is demonstrated how globally optimal bound values can be computed automatically using fractional programming theory. Two different solution approaches, the Charnes--Cooper transformation and Dinkelbach's algorithm are reviewed.
 Our implementation of the extended Davis--Kahan theorem is used to calculate bound values in three significant examples. First, a pairwise comparison is made of the spaces spanned by the eigenvectors of the graph shift operator matrices corresponding to different stochastic block model graphs. Second our bound is calculated on the distance of the spaces spanned by eigenvectors of the graph shift operators and their corresponding generating matrices in the stochastic blockmodel,  and, third,  on the sample and population covariance matrices in a spiked covariance model. Our extended bound values, using affine transformations, not only outperform the standard Davis--Kahan bounds in all examples where both theorems apply, but also demonstrate good performance in several cases where the standard Davis--Kahan theorem cannot be used. 
\end{abstract}
\begin{IEEEkeywords}
affine transform, Davis--Kahan theorem, comparing spaces spanned by eigenvectors, graph shift operator, fractional programs, PCA, stochastic blockmodel 
\end{IEEEkeywords}

\section{Introduction}
In the first part of this work (Paper I) we introduced the extended Davis--Kahan (DK) theorem for comparing two sets of consecutive and corresponding eigenvectors from any two symmetric matrices.

The extension incorporated a polynomial transform of one of the matrices which allows a relaxation and utilisation of the eigenvalue structure imposed by the standard DK theorem. As a result the bounds determined by the extended theorem are always at least as tight as those from the standard DK theorem. 

Paper I concentrated on the mathematics of the proposed approach.
In this second part of the work,  we turn our attention to computational issues, and also give some significant examples of applications of our extended DK theorem. The computational aspects are certainly challenging, and in this paper we only give a full discussion for the case of affine (linear) transformations; however, as exemplified by the applications,
the affine transformation can be very beneficial. 

We classify the problem of finding optimal affine transformation parameters for the extended DK bound as a fractional program. Fractional programing seems to be a less-well known class of optimization problems. 
For the history and recent advances in fractional programming see \cite{Frenk2008}.
An excellent overview of the solution and implementation of concave-convex fractional programs, the subclass of fractional programs which applies to our optimization problem,  is given in \cite{Isheden2012}.
In the wireless communication literature fractional programming approaches have recently found much use 
\cite{Cheung2016, DeMaio2011, Shen2018, Ye2018}.

Once we have chosen a solution method for our optimization problem we produce bounds in a range of different applications. 
The matrices considered are graph shift operators corresponding to graphs in a stochastic blockmodel and covariance matrices in a spiked covariance model, where the principal component analysis algorithm is well motivated. As stated in paper I, bounding the spaces spanned by the eigenvectors of the graph shift operators is particularly relevant to the signal processing community since they form the basis of the much utilised graph Fourier transform.

The remainder of this paper is structured as follows, in Section \ref{sec_affine_any_K} we summarise our extended DK theorem proved in Paper I using general polynomial transformations and then specialize to affine transformations. In Section \ref{sec_affine_optimization_problem} we present the problem of finding transformation parameters resulting in optimal bound values in our extended affine DK theorem. Furthermore, we prove that a trivial bound on the distance of the spaces spanned by the compared eigenvectors is always outperformed for comparisons of eigenvectors corresponding to either the largest or smallest eigenvalues. In Section \ref{sec_fractional_programming} we introduce concave-convex fractional programs, demonstrate that the problem from Section \ref{sec_affine_optimization_problem} can be brought into a fractional programming form and discuss two solution approaches, the Charnes--Cooper transformation and Dinkelbach's algorithm. Section \ref{sec_results} presents three significant examples in which our bound is computed. These examples are the comparison of the spaces spanned by the eigenvectors of the graph shift operators in \ref{sec_results_rep_mats}, of the graph shift operators to their corresponding generating matrices in the stochastic blockmodel in Section \ref{sec_results_rep_vs_B} and of the sample and population covariance matrices in a spiked covariance model in Section \ref{sec_results_PCA}. Our summary and conclusions are given in Section~\ref{sec:summary}.

\section{Problem Summary} \label{sec_affine_any_K}
\subsection{Bounds}
Let $\mathbb{V}_{n,\dimevec}$ denote the Stiefel manifold of $n \times r$ matrices with orthonormal columns.
Let $\mati, \matii \in \mathbb{R}^{n \times n}$ be symmetric matrices  with eigenvalues $\evali_1 \leq \evali_2 \leq \ldots \leq \evali_n$ and $\evalii_1 \leq \evalii_2 \leq \ldots \leq \evalii_n$ and corresponding eigenvectors $\{ w_1, w_2, \ldots, w_n\}$ and $\{ v_1, v_2, \ldots, v_n\}$, respectively. 
For $j\geq 0, r\geq 1,$ let the matrices holding the eigenvectors corresponding to $\dimevec$ consecutive eigenvalues of each matrix be denoted by $\evecsi_j = [w_{j+1}, \ldots, w_{j+\dimevec}]\in \mathbb{V}_{n,\dimevec}$ and $\evecsii_j = [v_{j+1}, \ldots, v_{j+\dimevec}]\in \mathbb{V}_{n,\dimevec}, $ the columns of which span the spaces $\mathcal{\evecsi}_j$ and $\mathcal{\evecsii}_j,$ respectively.
Under the stated conditions, it was shown in Paper I that
there exists a $Q \in O(\dimevec)$ such that
\begin{align} \label{eqn_cost_structure_transformed}
\|\evecsi_j - \evecsii_j Q \|_F & \leq c_{n,r}\left\| \evecsi_j \evecsi_j^T (I-\evecsii_j \evecsii_j^T)\right\|_2\nonumber\\
&  \leq  c_{n,r}  \frac{\left\|p(\mati) - \matii\right\|_2}{\delta_i},
\end{align}
where $O(r)$ is the group of $r\times r$ orthogonal matrices,   $c_{n,r}=\sqrt{2 \min(\dimevec, n-\dimevec)},$ $p(\cdot)$ is a polynomial matrix transformation and $\delta_i, i \in\{1,2\}$ are different values of the DK interval separation parameter corresponding to the two different DK interval choices. The usual or standard DK bounds follow by taking $p(\mati)=\mati$ in the second inequality. 
 In Paper I we saw that the first norm is directly related to the metric $\rho_1(\mathcal{\evecsi}_j, \mathcal{\evecsii}_j)= \inf_{R \in O(\dimevec)} \|\evecsi_j -\evecsii_j R \|_F,$ and the second relates to the metric $\rho_2(\mathcal{\evecsi}_j, \mathcal{\evecsii}_j) = \left\| \evecsi_j \evecsi_j^T \left(I - \evecsii_j \evecsii_j^T \right) \right\|_2.$

\subsection{Affine Transformations}
In this paper we focus on
affine matrix transformations $f(\mati) = c_1 \mati + c_0 I , \,\,(c_1, c_0 \in \mathbb{R})$ of one of the matrices under comparison (arbitrarily, $\mati$). 
So we set $p(\cdot)=f(\cdot)$
and for $j\geq 0, r\geq 1,$ from Paper I, the two DK interval choices  take the form
$S_1 = [a,b], S_2 = \mathbb{R}\backslash (a - \delta, b+\delta),$ with
\begin{align}
  a_1 \!&= \!\!\underset{i \in \{j+1, \ldots, j+\dimevec \}}{\min} \!f(\evali_i), b_1 \!=\!\! \underset{i \in \{j+1, \ldots, j+\dimevec \}}{\max} \!f(\evali_i), \nonumber\\
   \delta_1 &= \min\left(  \evalii_{j+\dimevec+1} - b_1, a_1 - \evalii_j\right) ;\label{eqn_interval_choice1} \\
 a_2 &= \evalii_{j+1},\, b_2 = \evalii_{j+\dimevec}\nonumber\\
 \delta_2 &\!=\! \min\!\!\left[  \underset{i \in \mathcal{A}_1}{\min}~ f(\evali_i)\! -\! b_2, a_2 \!-\!  \underset{i \in \mathcal{A}_2 }{\max}~ f(\evali_i) \right]. \label{eqn_interval_choice2} \end{align}
 The index sets ${\cal A}_1$ and ${\cal A}_2$ are given in Paper~I.

A main purpose of this paper is to solve the problem of optimizing the bound on the right-side of \eqnref{eqn_cost_structure_transformed} when $p(\cdot)=f(\cdot),$ i.e.,  minimize 
\begin{equation}\label{eq:thisisit}
c_{n,r}
\frac{\left\|f(\mati) - \matii\right\|_2}{\delta_i},\quad i\in \{1,2\},
\end{equation}
subject to the associated constraints given in Paper I. 
For an affine transformation, Constraints 2 of Paper I must be applied (Constraints 1 are subsumed). 
Constraints 2 take the form:
\begin{enumerate}[label = \Alph*]
\item \label{ass_poly_non_overlap_interval1} In the case of interval choice \eqnref{eqn_interval_choice1}, let the transformation parameters of $f(\cdot)$ be chosen such that, for given $j\geq 0, \dimevec\geq1$,
\begin{align*}
\delta_1&>0, \numberthis\label{eqn_delta1_pos}\\
a_1 - \evalii_{j+1} &< \delta_1,\numberthis\label{eqn_use_of_eigengap1} \\
\evalii_{j+\dimevec} - b_1 &< \delta_1.\numberthis\label{eqn_use_of_eigengap2} 
\end{align*}
\item \label{ass_poly_non_overlap_interval2} For interval choice \eqnref{eqn_interval_choice2}, let the transformation parameters of $f(\cdot)$ be chosen such that, for given $j\geq 0, \dimevec\geq1$,
 \begin{align}
\delta_2&>0, \label{eqn_delta2_pos}\\
 a_2 -  \min_{i \in \{j+1, \ldots, j+\dimevec\}} f(\evali_i) &< \delta_2, \label{eqn_overlap_left}\\
 \max_{i \in \{j+1, \ldots, j+\dimevec\}} f(\evali_i) - b_2 &< \delta_2.\label{eqn_overlap_right}
\end{align}
\end{enumerate}

Also as pointed out in Paper I, 
in addition to the two DK interval choices which give $\delta_1, \delta_2$ in 
\eqnref{eq:thisisit}, there are also the possibilities $c_1>0$ and $c_1<0$ in the affine transform.
So there are four different possible values for the DK interval separation $\delta,$ which are, 
\begin{align}
\!\!\!\delta_{1,+} &= \min\left(  \evalii_{j+\dimevec+1} - c_1\evali_{j+\dimevec}-c_0, c_1\evali_{j+1} +c_0 - \evalii_j\right)\label{eqn_delta_1+}\\
\!\!\!\delta_{1,-} &= \min\left(  \evalii_{j+\dimevec+1} - c_1\evali_{j+1}-c_0, 
c_1\evali_{j+\dimevec}+c_0 - \evalii_j\right)\label{eqn_delta_1-}\\
\!\!\!\delta_{2,+} &= \min\left(  c_1\evali_{j+\dimevec+1}+ c_0 - \evalii_{j+\dimevec}, \evalii_{j+1} - c_1\evali_j-c_0 \right)\label{eqn_delta_2+}\\
\!\!\!\delta_{2,-} &= \min\left(  c_1\evali_j+c_0 - \evalii_{j+\dimevec}, \evalii_{j+1} - c_1\evali_{j+\dimevec+1}-c_0 \right)\label{eqn_delta_2-}
\end{align}

\section{The bound as a numerical optimization problem} \label{sec_affine_optimization_problem}

In this section we frame the optimization of the affine bounds \eqnref{eq:thisisit}  as constrained optimization problems over the transformation parameters $c_1, c_0.$ Solving the optimization problems derived in this section results in the minimal bound  under affine transformations on the distance of the spaces spanned by two sets of eigenvectors.
Since there are four possibilities for the denominator in \eqnref{eq:thisisit}, namely
\eqnref{eqn_delta_1+}-\eqnref{eqn_delta_2-},
we have four optimization subproblems, 
which need to be solved in order to obtain the {\it overall} optimal DK bound. 

We study the solution of one in detail, and the rest follow analogously.
In \eqnref{eqn_affine_general_bound} we show the optimization subproblem for interval choice \eqnref{eqn_interval_choice1} for affine transformations with $c_1>0$, i.e., where $\delta = \delta_{1,+}$ in \eqnref{eqn_delta_1+}.
 Then $ a_1 \!= \!\!\underset{i \in \{j+1, \ldots, j+\dimevec \}}{\min} \!f(\evali_i)=f(\phi_{j+1})=c_1\phi_{j+1}+c_0$ and 
$b_1 \!=\!\! \underset{i \in \{j+1, \ldots, j+\dimevec \}}{\max} \!f(\evali_i)= f(\phi_{j+r})=c_1\phi_{j+r}+c_0.$ (Here $a_1=f(\phi_{j+1})$ and $b_1=f(\phi_{j+r})$ are given in Table I of Paper I, but follow from the preservation of ordering of the eigenvalues for an affine transform with $c_1>0$).
For interval choice \eqnref{eqn_interval_choice1} we need Constraints 2A, which we apply in \eqnref{eqn_affine_general_bound} via the last three rows of  the ``s.t''  (``subject to'') statement.
 From Corollary 1 of Paper I 
the objective function to be minimized is ${\left\|f(\mati) - \matii\right\|_2}/{\delta_{1,+}}. $
Hence the first subproblem takes the form
\begin{equation}\label{eqn_affine_general_bound}
\begin{aligned}
&\underset{c_1, c_0 }{\min}    
&& \frac{\left\| c_1\mati + c_0 I - \matii \right\|_2}{\delta_{1, +}},\\ 
& \text{s.t.}
&& c_1>0,\\
&&& \delta_{1,+} >0,\\
&&&\delta_{1,+}> \evalii_{j+\dimevec} -   c_1 \evali_{j+\dimevec} - c_0 ,\\
&&& \delta_{1,+}> c_1 \evali_{j+1} + c_0 - \evalii_{j+1}  .
\end{aligned}
\end{equation}

The remaining three subproblems follow a similar structure as \eqnref{eqn_affine_general_bound}, where $\delta$ equals \eqnref{eqn_delta_1-}, \eqnref{eqn_delta_2+} or \eqnref{eqn_delta_2-} and the values of the transformed spectrum are correspondingly taken from Table~I of Paper I.

\begin{remark} \label{rmk_omission_of_const_in_minimisation_problem}
In the objective function in \eqnref{eqn_affine_general_bound} we omit the constant $c_{n,r}$ present in \eqnref{eqn_cost_structure_transformed}, since it is inconsequential to the minimization. However, solutions of \eqnref{eqn_affine_general_bound} have to be multiplied by $c_{n,r}$ in order to obtain valid bounds.
\hfill $\lhd$
\end{remark}

In Proposition \ref{prop_trivial_upper_bound_and_limit_behaviour} we obtain a trivial upper bound on $\|\evecsi - \evecsii Q \|_F $ and demonstrate that, for comparisons of the first or the last $\dimevec$ eigenvectors, solutions of \eqnref{eqn_affine_general_bound} always approximate or improve upon this bound. Without considering the matrix transform, no such guarantees could be given.

\begin{proposition} \label{prop_trivial_upper_bound_and_limit_behaviour}
$\,$
\begin{enumerate}
\item Let $W, V\in \mathbb{V}_{n,\dimevec}.$ 
Then, for any $j\geq 0,$ there exists $Q \in O(\dimevec)$ such that,
\begin{equation} \label{eqn_trivial_bound}
\|\evecsi_j - \evecsii_j Q \|_F \leq c_{n,r}.
\end{equation}
\item  When comparing spaces spanned by the $r$ eigenvectors corresponding to either the $r$ largest or $r$ smallest eigenvalues,
the bound produced from  \eqnref{eqn_affine_general_bound} always approximates or improves upon \eqnref{eqn_trivial_bound}. 
\end{enumerate}

\end{proposition}

\begin{IEEEproof}
We begin by proving part 1 of the proposition.
From Lemma 2 of Paper I,
$$
\|\evecsi_j - \evecsii_j Q \|_F \leq c_{n,r}\left\| \evecsi_j \evecsi_j^T (I-\evecsii_j \evecsii_j^T)\right\|_2.
$$
Now, $\evecsi_j \evecsi_j^T$ and $(I-\evecsii_j \evecsii_j^T)$ are both projectors and hence all their eigenvalues are equal to either 0 or 1 \cite[p.~358]{Bernstein2009}. 
Therefore,
\begin{align*}
\|\evecsi_j- \evecsii_j Q \|_F &\leq c_{n,r}~\left\|\evecsi_j \evecsi_j^T \right\|_2 \left\| I-\evecsii_j \evecsii_j^T \right\|_2 \leq c_{n,r}.
\end{align*}
In Appendix~\ref{app_limit_of_bound} we show that, for comparisons of 
spaces spanned by the $r$ eigenvectors corresponding to either the $r$ largest or $r$ smallest eigenvalues,
 the objective function in \eqnref{eqn_affine_general_bound} approximates 1 as $c_0$ tends to either $\infty$ or $-\infty$.  
This result holds for all optimization problems of the form \eqnref{eqn_affine_general_bound}, where $j=0$ or $j = n-\dimevec$. 
As stated in Remark \ref{rmk_omission_of_const_in_minimisation_problem}, multiplying the objective function by $c_{n,r}$ yields valid bound values.
Hence, for $j \in \{0, n-\dimevec\}$ there exist large negative or positive values of $c_0$ such that the cost function \eqnref{eqn_affine_general_bound} results in bounds approximating \eqnref{eqn_trivial_bound}. 
For some problems of the form \eqnref{eqn_affine_general_bound}, an appropriate choice of transformation parameters may result in a bound smaller than 1 as shown in Section~\ref{sec_results}. 
Therefore, for $j \in \{0, n-\dimevec\}$ the bound produced from a solution of \eqnref{eqn_affine_general_bound} always approximates or improves upon \eqnref{eqn_trivial_bound}. 
\end{IEEEproof}

\begin{remark}
When considering other than the first or last $\dimevec$ eigenvectors, i.e., eigenvector comparisons with $1\leq j \leq n-1-r,$ we observe from Equations \eqnref{eqn_delta_1+}, \eqnref{eqn_delta_1-}, \eqnref{eqn_delta_2+} and \eqnref{eqn_delta_2-} that all four quantities  $\delta_{1,+}, \delta_{1,-}, \delta_{2,+}$ and $ \delta_{2,-}$ tend to $-\infty$ as $c_0 \rightarrow \pm \infty.$ Hence,  for eigenvector comparisons with $1\leq j \leq n-1-r,$ Constraints 2 are violated when $c_0 \rightarrow \pm \infty,$ i.e., $c_0 \rightarrow \pm \infty$ lies outside of the feasible set of \eqnref{eqn_affine_general_bound} and its related subproblems. Hence, for comparisons of spaces spanned by eigenvectors other than the first or last $r$, (ordered by corresponding eigenvalue magnitude), a similar statement to part 2 of Proposition \ref{prop_trivial_upper_bound_and_limit_behaviour} is not guaranteed.
\hfill $\lhd$
\end{remark}

\begin{remark}
The trivial upper bound in \eqnref{eqn_trivial_bound} applies for all values of $j.$ Suppose we are unable to find a set of affine transformation parameters such that $\|c_1 \mati + c_0 I - \matii\|_2<\delta.$ This means there exists no affine transformation reducing the distance of the two matrices  in the two norm to a value less than the distance of the relevant eigenvalues. In this case, the trivial bound in \eqnref{eqn_trivial_bound}  should be used to bound the distance of the spaces spanned by the $\dimevec$ consecutive eigenvectors of the two matrices instead of the bound resulting from the solutions of \eqnref{eqn_affine_general_bound} and its related subproblems.  
\hfill $\lhd$
\end{remark}

\section{Calculating the bound in practice} \label{sec_fractional_programming}

\subsection{Fractional Programming}
Here we show that the optimization subproblem \eqnref{eqn_affine_general_bound} and its related subproblems can be efficiently solved using fractional programming theory.

Ratio optimization problems are commonly called \textit{fractional programs} \cite{Frenk2008}. Hence, the optimization problem \eqnref{eqn_affine_general_bound} of choosing the affine transformation parameters $(c_1, c_0) \in \mathbb{R}^+ \times \mathbb{R}$ resulting in a minimal bound is a fractional program. 
We now describe the properties of fractional programs and their solutions. Then we transform \eqnref{eqn_affine_general_bound} to fit the standard class of concave-convex fractional programs and discuss the implementation of its solution.

Firstly, we formally define fractional programs.

\begin{definition} \label{defn_fractional_program}
\cite{Isheden2012} A general nonlinear fractional program has the form,
\begin{equation}\label{eqn_fractional_program_defn}
\begin{aligned}
&\underset{\mathbf{x}}{\max}    
&& \frac{g_1(\mathbf{x})}{g_2(\mathbf{x})},\,\,\text{s.t. } \mathbf{x} \in \mathcal{S},
\end{aligned}
\end{equation}
where $\mathcal{S} \subseteq \mathbb{R}^m$, $g_1, g_2: \mathcal{S}\rightarrow\mathbb{R}$ and $g_2(\mathbf{x})>0$. Problem \eqnref{eqn_fractional_program_defn} is called a concave-convex fractional program if $g_1$ is concave, $g_2$ is convex, and $\mathcal{S}$ is a convex set; additionally $g_1(\mathbf{x}) \geq 0$ for $\mathbf{x} \in \mathcal{S}$ is required, unless $g_2$ is affine.
\hfill $\lhd$
\end{definition}

In \cite{Frenk2008}, concave-convex fractional programs are referred to as concave fractional programs. An excellent overview of concave-convex fractional programs is given in \cite{Isheden2012}, with a focus on wireless communication.

\begin{remark} \label{rmk_fractional_optimization_global_optima}
For concave-convex fractional programs, a powerful and useful practical result is that  any local maximum is a global maximum \cite{Frenk2008}. 
\hfill $\lhd$
\end{remark}

When discussing the solution of concave-convex fractional programs the concept of equivalence of optimization problems is essential.

\begin{definition}
\cite[p.~130]{Boyd2004} define two optimization problems as \textit{equivalent} if the solution of one problem can be readily obtained given the solution of the other problem and vice versa.
\hfill $\lhd$
\end{definition}

Furthermore, we make use of the standard definition of the feasible set of an optimization problem.

\begin{definition}
\cite[p.~127]{Boyd2004} The \textit{feasible set} of an optimization problem is equal to the set of points which satisfy all the constraints of the optimization problem.
\hfill $\lhd$
\end{definition}

\subsection{Creating a Concave-convex Fractional Program}
Subproblem \eqnref{eqn_affine_general_bound} can be transformed to fall into the class of concave-convex fractional programs.
As pointed out in
\cite{Schaible1976}, 
\begin{equation} \label{eqn_equivalence_of_min_and_max}
\underset{\mathbf{x} \in \mathcal{S}}{\max} \left(\frac{g_1(\mathbf{x})}{g_2(\mathbf{x})}\right) = \frac{1}{\underset{\mathbf{x}\in \mathcal{S}}{\min}\left(\frac{g_2(\mathbf{x})}{g_1(\mathbf{x})}\right)}.
\end{equation}
Using \eqnref{eqn_equivalence_of_min_and_max} we find that solving \eqnref{eqn_affine_general_bound} is equivalent to solving,
\begin{equation}\label{eqn_affine_general_bound_concave_convex}
\begin{aligned}
&\underset{c_1, c_0 }{\max}    
&& \frac
{\delta_{1, +}}
{\left\| c_1 \mati + c_0 I - \matii \right\|_2}
,\\ 
& \text{s.t.}
&& c_1>0,\\
&&& \delta_{1,+} >0,\\
&&&\delta_{1,+}> \evalii_{j+\dimevec} -   c_1 \evali_{j+\dimevec} - c_0 ,\\
&&& \delta_{1,+}> c_1 \evali_{j+1} + c_0 - \evalii_{j+1}  .
\end{aligned}
\end{equation}

The bound value is found by transforming the solution of \eqnref{eqn_affine_general_bound_concave_convex} according to \eqnref{eqn_equivalence_of_min_and_max}.

We now demonstrate that the subproblems such as \eqnref{eqn_affine_general_bound_concave_convex}, are concave-convex fractional programs. 

\begin{lemma}
The $\delta$'s in Equations \eqnref{eqn_delta_1+}, \eqnref{eqn_delta_1-}, \eqnref{eqn_delta_2+} and \eqnref{eqn_delta_2-}  are all concave.
\end{lemma}
\begin{IEEEproof}
All the $\delta$'s in Equations \eqnref{eqn_delta_1+}--\eqnref{eqn_delta_2-}  are equal to the minimum of two affine functions of the transformation parameters. 

Let $f_1(x), f_2(x), x \in \mathbb{R}$ be affine functions. Then, $-f_1(x)$ and $-f_2(x)$ are still affine functions. Affine functions can be thought of as either convex or concave \cite[p.~67]{Boyd2004} and further the pointwise maximum of convex functions is convex \cite[p.~80]{Boyd2004}. Hence, $\max(-f_1(x), -f_2(x))$ is a convex function. 

If $f(x), x \in \mathbb{R}$ is a convex function, $-f(x)$ is concave \cite[p.~67]{Boyd2004}. Therefore, $-\max(-f_1(x), -f_2(x))$ is concave. But $-\max(-f_1(x), -f_2(x)) = \min(f_1(x), f_2(x)),$  is of the same form as the $\delta$'s, therefore $\delta$ is concave. 
\end{IEEEproof}

All four subproblems share the denominator $\left\| c_1 \mati + c_0 I - \matii \right\|_2$ which is easily shown to be convex via the triangle inequality. Furthermore, we require the denominator in \eqnref{eqn_affine_general_bound_concave_convex} to be strictly positive. 
Since $\left\| c_1 \mati + c_0 I - \matii \right\|_2$ is not affine we additionally require the numerator of \eqnref{eqn_affine_general_bound_concave_convex} to be positive on its feasible set; this is ensured by the 
$\delta$'s being positive on this set.  
Hence, \eqnref{eqn_affine_general_bound_concave_convex} and the remaining 3 subproblems are elements of the class of concave-convex fractional programs as in Definition \ref{defn_fractional_program}.

\subsection{Solving a Concave-convex Fractional Program}
Several general approaches to solving concave-convex fractional programs are presented in \cite{Isheden2012} and in \cite{Shen2018} fractional programming in the context of multiple-ratio problems is discussed.

In Section \ref{sec_fractional_programming_Charnes_Cooper}, we discuss the parameter-free approach where an equivalent convex problem is obtained through transformation of the optimization parameters; the transformation used is commonly referred to as the Charnes--Cooper transformation.  This transformed problem only needs to be solved once. 
In \cite{Cheung2016} the Charnes--Cooper transformation is used in the optimization of the energy spectral efficiency of a communication network and in \cite{DeMaio2011} the authors show that the maximum likelihood estimate of the steering direction of a signal for radar detection can be found by utilising the Charnes--Cooper transformation of a fractional programming problem.

In Section \ref{sec_fractioanl_programming_Dinkelbach}, we treat the parametric approach which introduces an additional parameter $\lambda$ to obtain an equivalent problem, which is not jointly convex in $(c_1, c_0)$ and $\lambda$. The equivalent problem, is however, convex in $(c_1, c_0)$ and monotone in $\lambda$. Therefore, we iteratively solve the convex problem for $(c_1, c_0)$ for a fixed $\lambda$ and update $\lambda$ using a Newton-Raphson step. This algorithm is credited to Dinkelbach \cite{Dinkelbach1967}. 
In \cite{Ye2018}  minimization of  the system outage probability in a communication network using Dinkelbach's algorithm is discussed, using a closed form solution to the problem at each iteration.

For computational reasons we mainly utilise the parameter-free approach, i.e., the Charnes--Cooper transformation. This follows advice  in \cite{Schaible1981} and \cite{Schaible1983}, who state that the iterative solution via Dinkelbach's algorithm is only to be preferred over the single Charnes--Cooper transformed problem, if the solution via Dinkelbach exploits the structure of the numerator and denominator of the fractional program which the Charnes--Cooper solution does not. For instance for quadratic fractional programs -- fractional programs with a quadratic numerator and denominator and affine constraints -- Dinkelbach's algorithm solves a quadratic program at every iteration, while the Charnes--Cooper transformation yields a concave problem. Therefore, if not many Dinkelbach iterations are necessary for convergence, then Dinkelbach's algorithm is to be preferred over the Charnes--Cooper approach for quadratic fractional programs. We find that for our problem \eqnref{eqn_affine_general_bound} both the Charnes--Cooper transformation and Dinkelbach's algorithm solve convex or concave problems. Therefore, we prefer the Charnes--Cooper solution method. However, most importantly, the results of the two different approaches agree in our simulations, as would be anticipated from 
\cite{Isheden2012} who showed that the optimality conditions of the two approaches are equivalent, so in theory  the results  should indeed not vary.

 \subsection{The Charnes--Cooper Transformation} \label{sec_fractional_programming_Charnes_Cooper}
Charnes and Cooper \cite{Charnes1962} proposed a variable transform for linear fractional programs -- fractional programs with an affine numerator and denominator and linear constraints. Schaible \cite{Schaible1974} generalised the transformation to concave-convex fractional programs.
\cite{Isheden2012} give a good recent summary of the transformation of concave-convex fractional problems. 
For the transformation of \eqnref{eqn_affine_general_bound_concave_convex}, appropriate transformation parameters are:
\begin{align} \label{eqn_Charnes_Cooper_parameters}
y_1 &= \frac{c_1}{\left\| c_1 \mati + c_0 I - \matii \right\|_2};\, y_2 = \frac{c_0}{\left\| c_1 \mati + c_0 I - \matii \right\|_2};\nonumber\\
&t = \frac{1}{\left\| c_1 \mati + c_0 I - \matii \right\|_2}.
\end{align}
Let
$$
\delta_{1,+}(t)  \!=\! \min\left( t \evalii_{j+\dimevec+1}\! -\! y_1 \evali_{j+r}\!-\!y_2,\right.
 \left. y_1\evali_{j+1} \!+\!y_2\!-\! t\evalii_j\right).
 $$
Transforming \eqnref{eqn_affine_general_bound_concave_convex} using the parameters in \eqnref{eqn_Charnes_Cooper_parameters} we obtain the following convex optimization problem:
\begin{equation}\label{eqn_affine_general_bound_linear}
\begin{aligned}
&\underset{y_1,y_2, t }{\max}  
&& \delta_{1,+}(t),\\ 
& \text{s.t.}
&& t>0,\\
&&& \left\| y_1 \mati + y_2 I - t \matii \right\|_2 \leq 1, \\
&&& y_1 > 0, \\
&&&  \delta_{1,+}(t) >0,\\
&&&  \delta_{1,+}(t)> t \evalii_{j+\dimevec} -   y_1 \evali_{j+\dimevec} - y_2,\\
&&& \delta_{1,+}(t)> y_1 \evali_{j+1} + y_2 - t \evalii_{j+1}.
\end{aligned} 
\end{equation}
In the original proposal of the transformation for linear fractional programs \cite{Charnes1962} the equality constraint $\left\| y_1 \mati + y_2 I - t \matii \right\|_2 = 1 $ was used. This constraint cannot be placed on concave-convex fractional problems, since convex optimization problems can only have linear equality constraints \cite[p.~191]{Boyd2004}. 
It is proved in 
\cite{Schaible1974} that for concave-convex fractional programs the constraints $\left\| y_1 \mati + y_2 I - t \matii \right\|_2 = 1 $ and $\left\| y_1 \mati + y_2 I - t \matii \right\|_2 \leq 1 $ are equivalent. Therefore, we work with the relaxed constraint $\left\| y_1 \mati + y_2 I - t \matii \right\|_2 \leq 1 $.

Note that \eqnref{eqn_affine_general_bound_linear} is not a linear program since the constraint $\left\| y_1 \mati + y_2 I - t \matii \right\|_2 \leq 1$ contains a non-linear function of the parameters. The constraint is however convex; therefore, \eqnref{eqn_affine_general_bound_linear} is a convex optimization problem.

We implement the 4 subproblems using the cvx package in MATLAB \cite{Grant2018, Grant2008}. cvx does not accept strict inequalities. 
We therefore solve relaxed subproblems, where the strict inequalities are relaxed to include their boundaries, and then check whether the obtained solutions satisfy the strict inequalities. We have found this approach  to work extremely well in practice with no convergence issues.

\subsection{Dinkelbach's Algorithm} \label{sec_fractioanl_programming_Dinkelbach}
Dinkelbach's algorithm was proposed in \cite{Dinkelbach1967}.
Equivalent to \eqnref{eqn_affine_general_bound_concave_convex} is the problem
\begin{equation}\label{eqn_affine_general_bound_Dinkelbach}
\begin{aligned}
&\underset{c,d }{\max}   
&& 
\delta_{1, +}
- \lambda \left\| c_1 \mati + c_0 I - \matii \right\|_2
,\\ 
& \text{s.t.}
&& c_1>0,\\
&&& \delta_{1,+} >0,\\
&&&\delta_{1,+}> \evalii_{j+\dimevec} -   c_1 \evali_{j+\dimevec} - c_0 ,\\
&&& \delta_{1,+}> c_1 \evali_{j+1} + c_0 - \evalii_{j+1}  .
\end{aligned}
\end{equation}
\cite{Dinkelbach1967} state that the algorithm can be initialised at a feasible point $(c_1, c_0)$, which is chosen such that the corresponding $\lambda = g_1(c_1, c_0)/g_2(c_1, c_0)$ is positive, or at $\lambda=0$. When we initialise at $\lambda=0$ then any feasible set of transformation parameters $(c_1, c_0)$ can be chosen for the initialisation. Therefore, we choose to always initialise at $\lambda=0$ and $(c_1, c_0)$ to be equal to their respective optima from the Charnes--Cooper algorithm in the corresponding subproblem.

As with the Charnes--Cooper implementation we utilise relaxed subproblems, where the strict inequality constraints are relaxed to include their boundaries. Then we check and report if any of the strict inequality constraints in the Dinkelbach implementation are violated.
We have found the solution of Dinkelbach's algorithm to agree with the solution of the Charnes--Cooper algorithm in all cases we tested.

Dinkelbach's scheme was extremely useful for checking that our implementation of the Charnes--Cooper scheme was correct, but it offered no advantages over the latter, and was much slower.

\section{Visualising the bound values:  three examples} \label{sec_results}
Problem~\eqnref{eqn_affine_general_bound} and its solution via 
 \eqnref{eqn_affine_general_bound_concave_convex}, (along with the three related optimization subproblems),
 can be used to calculate bounds on the distance of the spaces spanned by eigenvectors of \textit{any} two symmetric matrices satisfying 
 Assumption \assnonzerojtheigengaps\/ of Paper I. Therefore, we envisage that the affine transform could contribute tighter bounds in a range of fields where eigenvectors are used.
We highlight three such applications. 
In Section \ref{sec_results_rep_mats} we study our bound on the distance of the spaces spanned by the eigenvectors of the three graph shift operator matrices. Then, in Section \ref{sec_results_rep_vs_B} we will apply the bound to the comparison of the graph shift operator matrices to their respective generating matrices in the stochastic blockmodel. Our final example application 
in Section \ref{sec_results_PCA} is in a principal component analysis setting, where we compare the space spanned by the eigenvectors of the sample covariance matrix and its corresponding population covariance matrix in a spiked covariance model.

Throughout Sections \ref{sec_results_rep_mats} and \ref{sec_results_rep_vs_B} we will generate networks from the stochastic blockmodel, introduced next.

\subsection{The Stochastic Block Model}
The stochastic blockmodel, which is widely used in the networks literature \cite{Holland1983, KarrerNewman11, Lei2015}, allows us to encode a block structure in a random graph via different probabilities of edges within and between node-blocks. The definition and parametrisation below is adapted from \cite{Lei2015}.
\begin{definition}
Consider a graph with node set $\{v_1, \ldots, v_n\}.$ Split this node set into $K$ disjoint blocks denoted ${\mathcal{B}_1, \ldots, \mathcal{B}_K }.$ We encode block membership of the nodes via a membership matrix $M \in \{0,1\}^{n \times K}$, where $M_{i,j} = 1$ if $v_i \in \mathcal{B}_j$ and $M_{ij} = 0$ otherwise. Finally, we fix the probability of edges between blocks to be constant and collect these probabilities in a probability matrix $P \in [0,1]^{K \times K}$, i.e., for nodes $v_i \in \mathcal{B}_l$ and $v_j \in \mathcal{B}_m$ the probability of an edge between $v_i$ and $v_j$ is equal to $P_{l,m}.$
\end{definition}

Hence, the parameters of the stochastic blockmodel are $M \in \{0,1\}^{n \times K}$ and $P \in [0,1]^{K \times K},$ where the number of nodes $n \in \mathbb{N}$ and the number of clusters $K \in \mathbb{N}$ are implicitly defined via the dimensions of $M.$ We simulate graphs from this model by fixing these parameters and then sampling edges from Bernoulli trials. The Bernoulli parameter of the trial corresponding to the edge connecting $v_i$ to $v_j$ is given by entry $(i,j)$ of the matrix $B_A=M P M^T$. 

Since our results  apply only to symmetric matrices, we work with undirected graphs, which can be derived from a stochastic blockmodel by sampling the upper triangular half of the adjacency matrix from the stochastic blockmodel and then equating the lower triangular part of the adjacency matrix to the transpose of the upper triangle.

Throughout this section we take the matrix of edge probabilities to be composed of only two values. On the diagonal we have $p_{w}$ encoding the probability of edges within the different blocks to be the same for all blocks. Off-diagonal we have $p_{b}$ to encode the probability of edges between nodes in different blocks. For example, for $K=3,$ $P$ takes the form:
\begin{equation} \label{eqn_structure_of_P}
P = 
\left(\begin{smallmatrix} 
 p_{w} & p_{b} & p_{b}\\
 p_{b} & p_{w} & p_{b}\\
 p_{b} & p_{b} & p_{w}
\end{smallmatrix}\right).
\end{equation}

\subsection{Scaling}\label{subsec:scaling}
For affine transforms, from \eqnref{eqn_cost_structure_transformed} and \eqnref{eq:thisisit} the bound of interest is the right-hand-side of
\begin{align*}
\left\|\evecsi_j - \evecsii_j Q\right\|_F &\leq  c_{n,r} \left\|\evecsi_j \evecsi_j^T - (I- \evecsii_j \evecsii_j^T)\right\|_2  \\
&\leq c_{n,r} \frac{\left\|c_1 \mati +c_0 I - \matii\right\|_2}{\delta_i}, 
\end{align*}
where $\delta_i$ equals $\delta_{1,+}, \delta_{1,-}, \delta_{2,+}$ or $\delta_{2,-},$ from Equations \eqnref{eqn_delta_1+}-\eqnref{eqn_delta_2-}, which generate the four different optimization subproblems of the form \eqnref{eqn_affine_general_bound_concave_convex}. 

In this section we divide all bound values and attained values by $c_{n,r},$  i.e., we consider instead
\begin{align*}
c_{n,r}^{-1} \left\|\evecsi_j - \evecsii_j Q\right\|_F &\leq  \left\|\evecsi_j \evecsi_j^T - (I- \evecsii_j \evecsii_j^T)\right\|_2  \\
&\leq  \frac{\left\|c_1 \mati +c_0 I - \matii\right\|_2}{\delta_i}, \numberthis \label{eq:sharp}
\end{align*}
This rescaling has the advantage that, independent of $n$ and $\dimevec,$ our bound values are on the same relative scale. Furthermore, the trivial bound derived in Proposition \ref{prop_trivial_upper_bound_and_limit_behaviour}
 corresponds to the upper bound of 1 in all plots, rather than the value $c_{n,r},$ which would vary across the different simulations.
 
 In what follows, (scaled) attained  distance in the metric 
 $\rho_1(\mathcal{\evecsi}_j, \mathcal{\evecsii}_j)$ refers to the quantity $c_{n,r}^{-1}\inf_{R \in O(\dimevec)} \|\evecsi -\evecsii R \|_F$  and (scaled) attained distance in the metric $\rho_2(\mathcal{\evecsi}_j, \mathcal{\evecsii}_j)$
 refers to 
$\left\|\evecsi_j \evecsi_j^T - (I- \evecsii_j \evecsii_j^T)\right\|_2 .$
  (For the first of these we recall from Paper I that, when calculating distances, finding the matrix $Q$ for which the infimum is attained can be avoided by the use of canonical angles.)

\subsection{Different Pairs of Graph Shift Operator  Matrices} \label{sec_results_rep_mats} 

In this section we calculate the bound on the distance of the spaces spanned by the eigenvectors of the graph shift operator matrices. 
Recall that the largest eigenvalues of the adjacency matrix correspond to the smallest eigenvalues of the Laplacians. Hence, in two of the three presented  comparisons we will compare spaces spanned by  eigenvectors corresponding to eigenvalues on opposing ends of the eigenvalue spectrum. Therefore, in the majority of cases presented in this section the standard DK Theorem  does not apply, while bound values can be obtained via our extended DK Theorem.

 Throughout this section we consider stochastic blockmodels with $K=3,$ where every block is composed of equally many nodes. The parameters identifying the compared eigenvectors are $j=1$ and $r=2$. This choice of $j$ and $\dimevec$ is motivated by the fact that the first eigenvector of the Laplacian matrices is a constant vector and is therefore not informative in the recovery of the blocks in the stochastic blockmodel. 
 Hence, we are comparing spaces spanned by eigenvectors corresponding to the {\it second and third largest}\/ eigenvalue of the adjacency matrix to those corresponding to the {\it second and third smallest}\/ eigenvalues of the two graph Laplacians.

\begin{figure}[t!]
\begin{center}
\includegraphics[scale=0.6,clip]{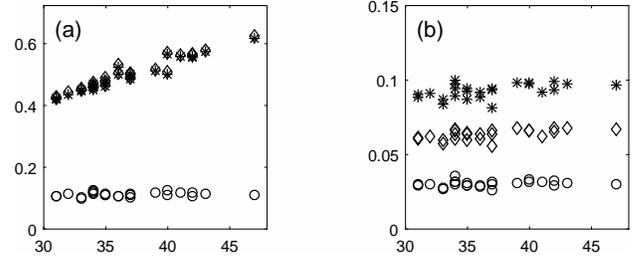}
\caption{
Different pairs of graph shift operator  matrices.
(a) Bound values plotted against the degree extreme differences ($x$-axis) of their corresponding network. The stars correspond to the comparison of spaces spanned by eigenvectors of $A$ and $L,$ diamonds correspond to the comparison of  $L$ and $\lsym$ and circles correspond to the comparison of  $A$ and $\lsym.$  
(b) The attained distances in the metric $\rho_1(\espacei_1, \espaceii_1)$ (see Paper I)
are plotted against the degree extreme differences ($x$-axis) of their corresponding network. All values have been rescaled as discussed in the text.}
\label{fig_sim1_degree_extreme_difference_vs_bd_and_attained}
\end{center}
\end{figure}

In Figs.~\ref{fig_sim1_degree_extreme_difference_vs_bd_and_attained} (a) and (b) we plot the bound and attained values arising from the comparison of the spaces spanned by eigenvectors of the three graph shift operator matrices against the degree extremes of the corresponding graphs. We simulated $n_{sim}=25$ stochastic blockmodels with equal parameters $n=300, K=3, (p_{b}, p_{w}) = (0.1,0.6).$

In Fig.~\ref{fig_sim1_degree_extreme_difference_vs_bd_and_attained}(a) we observe that the bound grows with a growing degree extreme difference for the comparisons of $L$ with $A$ and $\lsym$, while the bound remains relatively constant across different degree extremes when $A$ and $\lsym$ are compared.  The bound values for comparisons of $L$ with $A$ and $\lsym$ are very close, differing only by very small amounts. The bound arising from the comparison of $A$ and $\lsym$ attains much lower values than the bound values arising from the other two comparisons. Hence, using an affine transformation we are able to obtain a very small bound on the difference of spaces spanned by the eigenvectors of $A$ and $\lsym$, which suggests that they are very close. Not only are they close to each other, they also produce extremely similar bounds on the space spanned by the eigenvectors when individually compared to the eigenvectors of $L.$ 

In Fig.~\ref{fig_sim1_degree_extreme_difference_vs_bd_and_attained}(b) we 
observe that the attained distances of the three comparisons remain rather constant across the different degree extreme differences. 

The results in Fig.~\ref{fig_sim1_degree_extreme_difference_vs_bd_and_attained} show that for the comparison of $A$ and $\lsym$ an affine transformation is sufficient to remove the dependence of the bound value on the degree extreme difference. In contrast the much higher bounds for the comparison of the eigenvectors of $L$ with $A$ and $\lsym$ still depends on the degree extreme difference and the affine transformation was not sufficient to remove this dependence.

\begin{figure}[t!]
\begin{center}
\includegraphics[scale=0.53,clip]{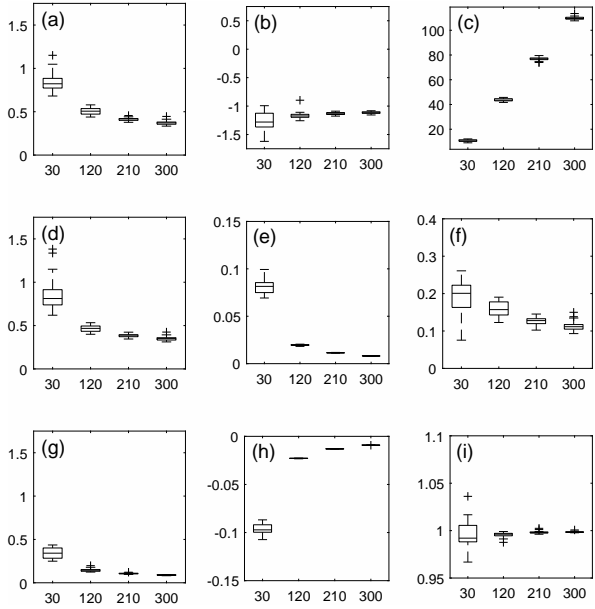}
\caption{
Different pairs of graph shift operator  matrices and
the effect of a growing number of nodes ($x$-axis).
First column: boxplots of the bound values of the three possible pairwise comparisons; second column: boxplots of optimal transformation parameter $c_1$; third column: boxplots of 
optimal transformation parameter $c_0$.
First row: comparison of $A$ and $L;$ second row:  $L$ and $\lsym;$ third row: $A$ and $\lsym.$  }
\label{fig_sim1_bound_vs_n_boxplot}
\end{center}
\end{figure}

In Fig.~\ref{fig_sim1_bound_vs_n_boxplot}, we observe the effect of a growing number of network nodes $n$ on our bound. For each value of $n \in \{30, 120, 210, 300\}$ we simulated $25$ stochastic blockmodels with equal parameters $K=3, (p_{b}, p_{w}) = (0.1,0.9).$
In all plots the four values of $n$ are displayed on the $x$-axis. The first column of plots in Fig.~\ref{fig_sim1_bound_vs_n_boxplot} displays boxplots of the bound values of the three possible pairwise comparisons, while the second and third columns display boxplots of the optimal affine transformation parameters $c_1$ and $c_0,$ respectively. 
The rows of plots show  the comparison of spaces spanned by eigenvectors of $A$ and $L,$ (first),  $L$ and $\lsym,$ (second)  and $A$ and $\lsym,$ (third).  
Just the 25 samples from each stochastic blockmodel parametrisation 
were sufficient to reveal general trends in the 9 plots. 

From Fig.~\ref{fig_sim1_bound_vs_n_boxplot} (a), (d) and (g) we immediately see that the bound values decrease as the number of nodes in the stochastic blockmodels grows indicating that as $n$ grows the differences between the spaces spanned by the eigenvectors of the graph shift operator matrices decrease. 
From Fig.~\ref{fig_sim1_bound_vs_n_boxplot} (b), (c), (e), (f), (h) and (i) we observe the majority of optimal affine transformation parameters depend on $n.$ 
For the comparison of  $A$ and $L,$ in Fig.~\ref{fig_sim1_bound_vs_n_boxplot}(c), we find the optimal additive parameter $c_0$ to grow from roughly 10 for $n=30$ to roughly 110 for $n=300,$ while the multiplicative parameter  in Fig.~\ref{fig_sim1_bound_vs_n_boxplot}(b) remains almost constant for all values of $n$. Comparison of the Laplacian matrices $L$ and $\lsym,$ shows both transformation parameters vary slightly with $n,$  (Figs.~\ref{fig_sim1_bound_vs_n_boxplot}(e) and (f)). Finally, comparison of  $A$ and $\lsym,$ shows the additive parameter $c_0$ to remain mostly constant with changing $n$  (Fig.~\ref{fig_sim1_bound_vs_n_boxplot}(i)) while the multiplicative parameter $c_1$ grows with $n$ (Fig.~\ref{fig_sim1_bound_vs_n_boxplot}(h)). 
The magnitude by which the transformation parameters change as $n$ grows is 
clearly quite variable.

We were unable to run the simulation displayed in Fig.~\ref{fig_sim1_bound_vs_n_boxplot} beyond $n=300$ within reasonable computation time. 
The  times for the simulations with the four different values of $n$ were,  $3, 21, 118$ and $ 472$ minutes, respectively. For each value of $n,$ 300 convex optimization problems were solved since each bound on the three possible graph shift operator matrix comparisons was calculated for 25 different stochastic blockmodels per value of $n$ and  calculation of each bound involves the solution of 4 different convex optimization subproblems.  

\subsection{Graph Shift Operator Matrices and Generating Matrices} \label{sec_results_rep_vs_B}

In this section we compare (i) spaces spanned by eigenvectors of the graph shift operator  matrices to (ii) spaces spanned by eigenvectors of their corresponding generating matrices in the stochastic blockmodel $B_A = M P M^T - \mathrm{diag}(M P M^T), B_L = \mathrm{diag}(B_A \mathbf{1}_n) - B_A$ and $B_{\lsym} = \mathrm{diag}(B_A \mathbf{1}_n)^{-1/2} B_L \mathrm{diag}(B_A \mathbf{1}_n)^{-1/2}.$ (Here $\mathbf{1}_n$ is a column vector of ones with $n$ entries and the term $-\mathrm{diag}(M P M^T)$ in the calculation of $B_A$ ensures that our stochastic blockmodels do not have self-loops.) 
It is natural to compare the spaces spanned by these eigenvectors, since  consistency and rate of convergence of different methods, based on the eigenvectors of the graph shift operator matrices in a stochastic blockmodel setting, can be demonstrated \cite{Cape2019, Eldridge2018, Rohe2011}.

\begin{figure}[t!]
\begin{center}
\includegraphics[scale=0.55,clip]
{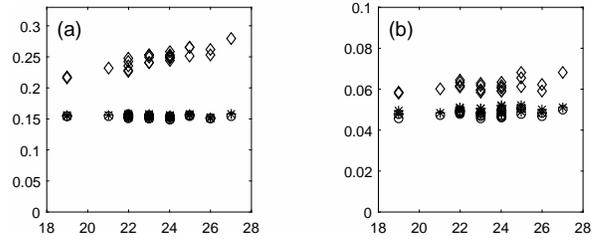}
\caption{
Comparison of the spaces spanned by eigenvectors of the graph shift operator matrices to those of their respective generating matrices, plotted against the degree extremes ($x$-axis) of their corresponding network .
The stars correspond to the comparison of eigenvectors of $A$ to $B_A,$ diamonds correspond to the comparison of eigenvectors of $L$ and $B_L$ and circles correspond to the bound values arising from the comparison of $\lsym$ and $B_{\lsym}.$  (a) bound values,  
and (b) attained distances in the metric $\rho_1(\espacei_0, \espaceii_0),$ (see Paper I).
All displayed values are rescaled as discussed in the text.
}
\label{fig_degree_extreme_difference_vs_bound_and_attained_repvsB}
\end{center}
\end{figure}

In Figs.~\ref{fig_degree_extreme_difference_vs_bound_and_attained_repvsB}(a) and (b) the bound and attained values from the comparison of the spaces spanned by the eigenvectors of the three graph shift operator matrices to the eigenvectors of their respective generating matrices are plotted against the degree extremes of the corresponding graphs. 
In the adjacency matrix comparison, the eigenvectors corresponding to the three largest eigenvalues are compared, while for the Laplacians we concern ourselves with the eigenvectors corresponding to the three smallest eigenvalues.
$25$ realisations of a stochastic blockmodel with parameters $n=210, K=3, (p_{b}, p_{w}) = (0.1,0.9)$ were simulated. In this comparison we included the first eigenvector of the matrices under comparison, i.e., we chose $j=0$ and $\dimevec =3.$

We see in Fig.~\ref{fig_degree_extreme_difference_vs_bound_and_attained_repvsB} that the ordering in magnitude of the attained values is reflected in the bound values, with the bounds on the unnormalised Laplacian comparison $(L, B_L)$ taking the largest values. 
We see that both the attained and the bound values in the comparison of the spaces spanned by eigenvectors of $L$ and $B_L$ seem to grow with growing degree extreme difference, which is not the case for the comparisons involving $A$ and $\lsym.$ 

\begin{figure}[t!]
\begin{center}
\includegraphics[scale=0.53,clip]
{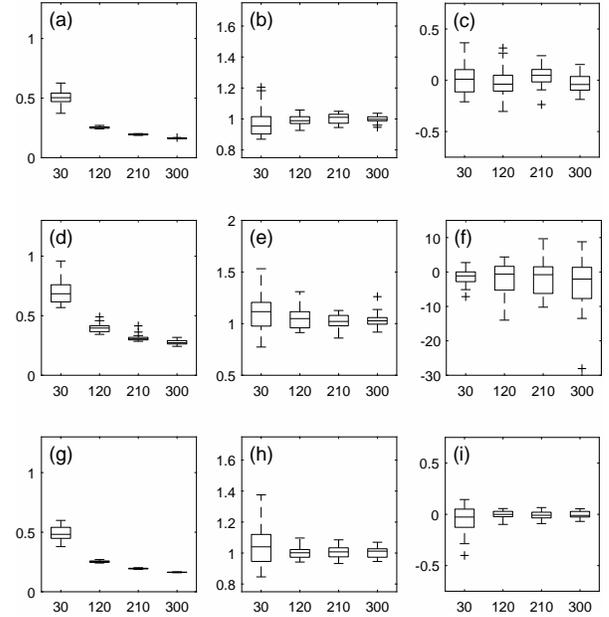}
\caption{
Graph shift operator matrices and generating matrices:
the effect of a growing number of nodes ($x$-axis).
First column: boxplots of the bound values; second column: boxplots of optimal transformation parameter $c_1$; third column: boxplots of 
optimal transformation parameter $c_0$.
The first row is for the comparison of spaces spanned by eigenvectors of $A$ and $B_A,$ the second row, $L$ and $B_L,$ and the third row, $\lsym$ and $B_{\lsym}.$}
\label{fig_bound_vs_n_boxplot_repvsB}
\end{center}
\end{figure}

In Fig.~\ref{fig_bound_vs_n_boxplot_repvsB} 
we observe the effect of a growing number of network nodes $n$ on our bound. For each value of $n \in \{30, 120, 210, 300\}$ we simulated $25$ stochastic blockmodels with  parameters $K=3, (p_{b}, p_{w}) = (0.1,0.8).$
The transformation parameters are roughly centred around the parameters of the identity transformation, $f(x)=x,$ i.e., $c_1=1$ and $c_0=0$. 
Interestingly, the variance of the additive parameter $c_0$ seems to be increasing with increasing $n$ for the comparison $(L,B_L).$ For all other displayed values in Fig.~\ref{fig_bound_vs_n_boxplot_repvsB}  we find the variance of the observed values to decrease as the number of nodes in the network grows. 
 As the number of nodes, $n,$ in the stochastic blockmodel grows, the bound values of all three comparisons decrease.

\begin{figure}[t!]
\begin{center}
\includegraphics[scale=0.53,clip]{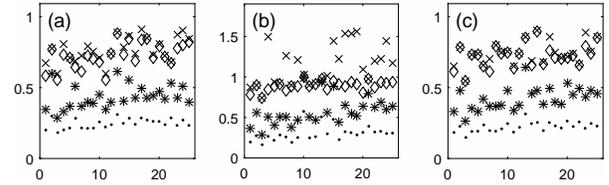}
\caption{
Graph shift operator matrices and generating matrices:
comparison of spaces spanned by eigenvectors of (a) $A$ and $B_A$, (b)  $L$ and $B_L$,  (c)  $\lsym$ and  $B_{\lsym}.$ 
 The dots represent the attained  distances in the metric $\rho_1(\espacei_0, \espaceii_0)$ 
 and the stars represent attained distances in the metric $\rho_2(\espacei_0, \espaceii_0)$ 
 (see Paper I).
 The diamonds show
 our sharpened DK bound using an affine  matrix transformation, and the x's represent usual DK bound values. 
All displayed values were rescaled as discussed in the text.}
\label{fig_attained_vs_bound_repvsB}
\end{center}
\end{figure}

For the spaces spanned by the leading eigenvectors of the graph shift operator matrices compared to their corresponding generating matrices, 
Fig.~\ref{fig_attained_vs_bound_repvsB} shows the usual DK bounds 
(Theorem~3 of Paper I), our sharpened bounds (Theorem~5 of Paper I and \eqnref{eq:sharp}) 
and the attained values; all were standarized.
 Here $25$ realisations of a stochastic blockmodel with parameters  $n = 30, K=3, (p_{b}, p_{w}) = (0.1,0.6)$ were generated. In all three comparisons, our sharpened bound values improve on the usual DK bound values. In the case of the unnormalised Laplacian several of  the usual DK bound values  are  greater than 1, therefore, even the trivial bound value of 1 (see Proposition~\ref{prop_trivial_upper_bound_and_limit_behaviour}
 and Section~\ref{subsec:scaling}) is tighter than the usual DK bound. In contrast, our bound produces values consistently lower than 1. 
 
 In addition to the attained distances of the spaces spanned by the eigenvectors in the metric $\rho_1(\espacei_0, \espaceii_0),$ we have shown the distance in the metric $\rho_2(\espacei_0, \espaceii_0)$ in Fig.~\ref{fig_attained_vs_bound_repvsB}, which is discussed in Paper I. As expected from Lemma 2 and Theorem 5 in paper I, we find the attained values $\rho_2(\espacei_0, \espaceii_0)$ to fall between the distances in the metric $\rho_1(\espacei_0, \espaceii_0)$ and our sharpened DK values. Of particular interest are the values attained in simulation number 10 in Fig.~\ref{fig_attained_vs_bound_repvsB}(b), where we find the distance of the spaces spanned by the eigenvectors to come very close to 1 in the metric $\rho_2(\espacei_0, \espaceii_0)$ and our sharpened bound to be very close to tight in this instance.

\subsection{Sample Covariance  and Population Covariance Matrices}\label{sec_results_PCA}

Our final example of the application of our sharpened DK bound is in the setting of 
Principal Component Analysis (PCA).
Consider $N \in \mathbb{N}$  independent, identically distributed samples $\{X_i\}$ from a multivariate normal distribution of dimension $p \in \mathbb{N},$ with mean $\mathbf{0}$ and covariance matrix $\Sigma.$
Then let $X$ be the $p\times N$ matrix with columns $X_i$ with $i \in \{1, \ldots, N\}$. We denote the sample covariance matrix by $\hat{\Sigma} = XX^T/N.$
When a low dimensional structure truly generates the covariance matrix, PCA is the correct tool to recover this low dimensional space. The standard PCA algorithm maps the data into the space spanned by the $\dimevec$ eigenvectors corresponding to the largest eigenvalues of $\hat{\Sigma}.$ In this setting it is of interest to study the convergence of the space spanned by these leading $\dimevec$ eigenvectors of $\hat{\Sigma}$ to the leading eigenvectors spanning the true low dimensional covariance space of $\Sigma$. 
Using a so-called spiked covariance model, we will apply our bound to the spaces spanned by eigenvectors corresponding to the largest eigenvalues of $\Sigma$ and $\hat{\Sigma}$.

The spiked covariance model was first introduced by \cite{Johnstone2001}, who described a phenomenon in real world data where the largest eigenvalues of the covariance matrix are separated by a large eigengap from the rest of the spectrum. \cite{Jung2009} proved consistency of the first $\dimevec$ eigenvectors of the estimated sample covariance towards the population covariance in the case of zero mean normally-distributed data and under certain conditions on the growth of the largest eigenvalues with growing dimensions $p$ and $N.$ 
The spiked covariance model has been found to be implied by the factor model \cite{Fan2013,ShenShen2016,Wang2017}, which models a multivariate time series as being driven by a few main factors. The factor model and consequently the spiked covariance model, find application to financial data. 

In our parametrisation of the spiked covariance model $r \in \mathbb{N}$ determines the dimension of the low-dimensional latent space in which $\Sigma$ is generated, $M \in \{0,1\}^{p \times \dimevec}$ encodes the latent dimension membership (similar to the stochastic blockmodel) and $P\in \mathbb{R}^{\dimevec \times \dimevec}$ encodes the correlations between latent dimensions.  $\Sigma$ has to be a valid covariance matrix, so must be  symmetric positive definite. As for the stochastic blockmodels, we chose $P$ to consist of only two values and take the form given in \eqnref{eqn_structure_of_P}.
We define the population covariance matrix as,
$$
\Sigma = M P M^T +  I.
$$
By building our covariance matrix like this, we  get a covariance matrix following the spiked covariance model, where $\dimevec$ eigenvalues are significantly larger than the rest of the spectrum (\cite{Johnstone2001, Jung2009}), the latter consisting of eigenvalues all equal to 1, as a result of adding $I$ into the covariance structure.
\begin{figure}[t!]
\begin{center}
\includegraphics[scale=0.53,clip]
{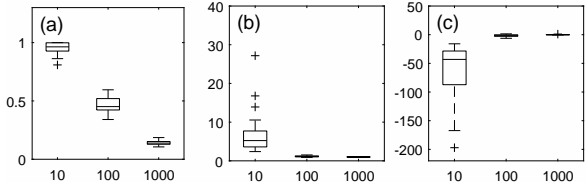}
\caption{
Effect of increasing sample size $N$  ($x$-axis) for the space spanned by eigenvectors of the sample covariance matrix versus the population covariance matrix.
(a) boxplots of bound values, (b) boxplots of optimal transformation parameter $c_1,$
(c) boxplots of optimal transformation parameter $c_0$.
 }
\label{fig_PCA_bound_vs_T}
\end{center}
\end{figure}

In Fig.~\ref{fig_PCA_bound_vs_T}(a) we plot our sharpened bound  
for spaces spanned by the eigenvectors corresponding to the largest 3 eigenvalues of $\Sigma$ and $\hat{\Sigma}$ and in Fig.~\ref{fig_PCA_bound_vs_T}(b) and (c) the corresponding optimal transformation parameters $c_1$ and $c_0$, respectively. 
Here $N \in \{10, 100, 1000\}, p = 60, (p_{b}, p_{w}) = (0.2, 0.8), j=0, \dimevec = 3.$ 
In Fig.~\ref{fig_PCA_bound_vs_T}(a) we see the bound decreases as the sample size $N$  grows. This makes sense: more samples should improve the estimation performance of the sample covariance matrix and therefore the distance of the subspaces spanned by the first three eigenvectors of the sample and population covariance matrices should decrease.  For $N=10$ samples we find that in a few cases the bound value 1 is attained. This situation was theoretically discussed in Proposition \ref{prop_trivial_upper_bound_and_limit_behaviour}. It is nice to see that we do indeed find the bound to converge to 1 with diverging transformation parameters in the worst case in practice. In plots (b) and (c) we find the transformation parameters to converge to the identity transformation as $N$ grows.

\begin{figure}[t!]
\begin{center}
\includegraphics[scale=0.52,clip]
{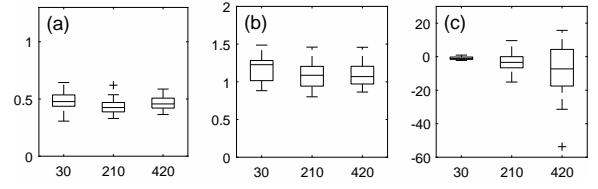}
\caption{Effect of increasing dimension $p$  ($x$-axis) for the space spanned by the eigenvectors of the sample covariance matrix versus population covariance matrix.
(a) boxplots of bound values, (b) boxplots of optimal transformation parameter $c_1,$
(c) boxplots of optimal transformation parameter $c_0$.
}
\label{fig_PCA_bound_vs_p}
\end{center}
\end{figure}

The format of Fig.~\ref{fig_PCA_bound_vs_p} follows Fig.~\ref{fig_PCA_bound_vs_T} with the difference that we keep $N$ fixed at 100 and study the behaviour of our bound as $p$ grows, $p \in \{30, 210, 420\}.$ Interestingly, our bound remains fairly constant for the values of $p$ considered. The transformation parameters hover around the identity transformation $f(x) =x$ with the uncertainty in the additive parameter $c_0$ increasing as $p$ increases.

\begin{figure}[t!]
\begin{center}
\includegraphics[scale=0.52,clip]
{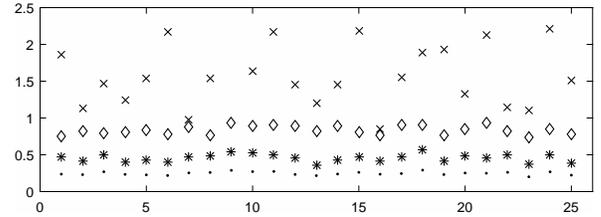}
\caption{Spaces spanned by first three eigenvectors of the sample covariance matrix versus the population covariance matrix. The dots represent attained distances in the metric $\rho_1(\espacei_0, \espaceii_0)$ and the stars represent attained distances in the metric  
$\rho_2(\espacei_0, \espaceii_0)$ (see Paper I) The diamonds show our sharpened DK bound using an affine matrix transformation. The x's show the usual DK bounds. All displayed values are rescaled as discussed in the text.}
\label{fig_PCA_atained_bound_Bhatia}
\end{center}
\end{figure}

Using  $N = 100, p = 60, (p_{b}, p_{w}) = (0.4, 0.6), j=0, r = 3,$ and 25 simulations,
Fig.~\ref{fig_PCA_atained_bound_Bhatia} 
shows that our extended or sharpened DK bounds 
improve upon the usual DK bounds by roughly a factor of 2. In Proposition \ref{prop_trivial_upper_bound_and_limit_behaviour} it was discussed that bound values above 1 are non-informative. In Fig.~\ref{fig_PCA_atained_bound_Bhatia} we observe our bound to consistently fall below 1, improving on the trivial bound; however, the usual DK bound  attains values consistently above 1.

\section{Summary and Conclusions}\label{sec:summary}

Paper I discussed the theory of polynomial transformations of $\mati, $ $p(\mati) = c_l \mati^l + c_{l-1} \mati^{l-1} + \ldots +c_1 \mati + c_0 I.$ Here we have concentrated on affine transformations which
 have the advantage of being monotone and hence the largest and smallest transformed eigenvalues in eigenvalue intervals are easily determined, e.g., 
$$
\underset{i \in \{1, \ldots, \dimevec\}}{\min} \left(f(\evali_i)\right) =  \begin{cases} 
      c_1 \evali_1 + c_0 & \text{for} \quad c_1 \geq 0, \\
      c_1 \evali_\dimevec + c_0 & \text{for} \quad  c_1 < 0.
   \end{cases}
$$
The overall minimal bound is thus found by considering 2 different optimization problems, where affine transforms with $c_1>0$ and $c_1<0$ are treated separately.

By way of contrast, 
higher order transformations are not monotone in general.
As a result, it is difficult to classify the resulting optimization problem to fall within a certain class of solvable optimization problems.
Hence, we have here restricted ourselves to the practical implementation of affine transformations. We found that in comparisons amongst the graph shift operator matrices, and of graph shift operator matrices with their generating matrices (in a stochastic blockmodel setting), our fractional programming implementation of the affine DK bounds
is superior to the standard DK bounds. The same was found when working with the eigenvectors of the sample and population covariance matrices in a spiked covariance model (for which  the PCA algorithm is a well motivated analysis tool).

\section*{Acknowledgment}
The work of Johannes Lutzeyer is supported by the EPSRC (UK).

\appendix{}

\subsection{Proof of  Part 2 of Proposition~\ref{prop_trivial_upper_bound_and_limit_behaviour}}\label{app_limit_of_bound}

In this Appendix we first work out the limit of bound \eqnref{eqn_affine_general_bound} as either $c_0 \rightarrow \infty$ or $c_0 \rightarrow -\infty$ for comparisons with $j=0$. Then we draw parallels to the comparisons involving the last $\dimevec$ eigenvectors, i.e., the case where $j=n-\dimevec.$

 We begin by producing a lower and an upper bound for the cost function in \eqnref{eqn_affine_general_bound}. Using the matrix triangle  and reverse-triangle inequalities
 $
 \big | ||A||-||B|| \big | \leq ||A-B|| $ and $ ||A+B||\leq ||A||+||B||
 $
 on the numerator of the cost function, gives
\begin{align}\label{eqn_app_upper_and_lower_bound}
\frac{\big | \left|c_0\right| \left\| I\right\|_2 - \left\| - c_1 \mati + \matii \right\|_2 \big |}{\delta} 
&\leq \frac{\left\| c_1 \mati + c_0 I - \matii \right\|_2}{\delta}\\
&\leq \frac{|c_0| \left\| I\right\|_2 + \left\| c_1 \mati - \matii \right\|_2}{\delta}.
\end{align}
But $\left\| I\right\|_2 =1$ so henceforth this term will be omitted.
Here $\delta$ is any of $\delta_{1,+}, \delta_{1,-}, \delta_{2,+}$ and $\delta_{2,-}$ corresponding to the four  subproblems which need to be solved in the affine case. 
For $j=0,$ we encounter the issue that $\evali_0$ and $\evalii_0$ are undefined and therefore, as is done in \cite[p.~317]{Yu2015} we set them equal to $-\infty.$ Then, 
\eqnref{eqn_delta_1+}, \eqnref{eqn_delta_1-}, \eqnref{eqn_delta_2+} and \eqnref{eqn_delta_2-} take the following form, 
\begin{align}
\delta_{1,+} &= \evalii_{\dimevec+1} -  c_1 \evali_{\dimevec} - c_0, \label{eqn_app_delta_1+} \\
\delta_{1,-} &= \evalii_{\dimevec+1} -  c_1 \evali_{1} - c_0, \label{eqn_app_delta_1-}\\
\delta_{2,+} &=  c_1 \evali_{\dimevec+1} + c_0 - \evalii_{\dimevec}, \label{eqn_app_delta_2+}\\
\delta_{2,-} &=  \evalii_{1}  - c_1 \evali_{\dimevec+1} - c_0. \label{eqn_app_delta_2-}
\end{align}

In order to work out the limit of the lower and upper bound in \eqnref{eqn_app_upper_and_lower_bound} as either $c_0 \rightarrow \infty$ or $c_0 \rightarrow -\infty$ we use l'Hopital's rule. For l'Hopital's rule to apply we require both the numerator and the denominator of the lower and upper bound in \eqnref{eqn_app_upper_and_lower_bound} to tend to 0 or $\infty$ as either $c_0 \rightarrow \infty$ or $c_0 \rightarrow -\infty$. We find that both numerators and their common denominators do indeed tend to $\infty$ as either $c_0 \rightarrow \infty$ or $c_0 \rightarrow -\infty$, i.e.,
\begin{align*}
\underset{c_0 \rightarrow \pm \infty }{\lim} \big| \left|c_0\right|   - \left\| - c_1 \mati + \matii \right\|_2 \big | &= \infty,\\
\underset{c_0 \rightarrow \pm \infty }{\lim} |c_0|   + \left\| c_1 \mati - \matii \right\|_2 &= \infty,\\
\underset{c_0 \rightarrow - \infty }{\lim} \delta_{1,+} = 
\underset{c_0 \rightarrow - \infty }{\lim} \delta_{1,-}= 
\underset{c_0 \rightarrow  \infty }{\lim} \delta_{2,+} &=
\underset{c_0 \rightarrow - \infty }{\lim} \delta_{2,-} = \infty.
\end{align*}

We start by evaluating the derivative with respect to $c_0$ of the numerator of both the lower and upper bound in \eqnref{eqn_app_upper_and_lower_bound} and of their denominators. 

Consider firstly $
({\partial}/{\partial c_0} )\left(\big | \left|c_0\right|   - \left\|  c_1 \mati - \matii \right\|_2 \big |\right) = \mathrm{sign}\left( \left|c_0\right|   - \left\|  c_1 \mati - \matii \right\|_2 \right) \mathrm{sign}(c_0) .
$
This derivative is:
$1$ for $ \left\|  c_1 \mati - \matii \right\|_2 < c_0;$
$-1$ for $0 < c_0 <  \left\|  c_1 \mati - \matii \right\|_2;$
$1$ for $-\left\|  c_1 \mati - \matii \right\|_2 < c_0 < 0,$ and
$-1$ for $c_0<- \left\|  c_1 \mati - \matii \right\|_2.$

Next we note
$({\partial}/{\partial c_0} )\left(|c_0|   + \left\| c_1 \mati - \matii \right\|_2\right) = \mathrm{sign}(c_0) =1$ for $c_0>0$ and $-1$ for $c_0<0.$

The derivative of the terms \eqnref{eqn_app_delta_1+}, \eqnref{eqn_app_delta_1-}, \eqnref{eqn_app_delta_2+} and \eqnref{eqn_app_delta_2-} follows trivially since they are linear functions of $c_0:$ 
\begin{equation*}
\frac{\partial}{\partial c_0}\delta_{1,+} = 
\frac{\partial}{\partial c_0}\delta_{1,-} = 
\frac{\partial}{\partial c_0}\delta_{2,-} = - 1;\quad
\frac{\partial}{\partial c_0}\delta_{2,+} = + 1.
\end{equation*}

For $\delta_{2,+},$ the limits of both the lower and the upper bound in \eqnref{eqn_app_upper_and_lower_bound} follow from l'H{\^ o}pital's rule:
\begin{align*}
\underset{c_0 \rightarrow   \infty }{\lim} \frac{\frac{\partial}{\partial c_0} \big | \left|c_0\right|   - \left\| - c_1 \mati + \matii \right\|_2 \big |}{\frac{\partial}{\partial c_0} \delta_{2,+}} &= 1\\
\qquad \xRightarrow{\text{l'H{\^ o}pital}}  \qquad  
\underset{c_0 \rightarrow   \infty }{\lim} \frac{ \big | \left|c_0\right|   - \left\| - c_1 \mati + \matii \right\|_2 \big |}{ \delta_{2,+}} &= 1;\\
\underset{c_0 \rightarrow   \infty }{\lim} \frac{\frac{\partial}{\partial c_0} \left(|c_0|   + \left\| c_1 \mati - \matii \right\|_2\right)}{\frac{\partial}{\partial c_0} \delta_{2,+}} &= 1\\
 \qquad \xRightarrow{\text{l'H{\^ o}pital}}   \qquad
\underset{d \rightarrow   \infty }{\lim} \frac{ \left(|d|   + \left\| c \mati - \matii \right\|_2\right)}{ \delta_{2,+}} &= 1.
\end{align*}
Since both the lower and the upper bound on \eqnref{eqn_affine_general_bound} tend to 1 as $c_0 \rightarrow   \infty,$  the sandwich theorem says  that 
$$
\underset{c_0\rightarrow  \infty }{\lim} \frac{ \left\| c_1 \mati + c_0 I - \matii \right\|_2}{ \delta_{2,+}} = 1.
$$
We also use l'H{\^ o}pital's rule and then the sandwich theorem to establish that for $\delta_{1,+},$ 
$$\underset{c_0 \rightarrow - \infty }{\lim} \frac{ \left\| c_1 \mati + c_0 I - \matii \right\|_2}{ \delta_{1,+}} = 1,
$$
and the same holds for $ \delta_{1,-}$ and $\delta_{2,-}.$

For $j=n-\dimevec,$ 
we set $\evali_{n+1}$ and $\evalii_{n+1}$ to equal $\infty,$ as done in \cite[p.~317]{Yu2015}, and then quantities \eqnref{eqn_delta_1+}, \eqnref{eqn_delta_1-}, \eqnref{eqn_delta_2+} and \eqnref{eqn_delta_2-} corresponding to the four  subproblems in the affine case are, 
\begin{align*}
\delta_{1,+} &= c_1 \evali_{n-\dimevec+1} + c_0 - \evalii_{n-\dimevec}, \\
\delta_{1,-} &= c_1 \evali_{n} + c_0 - \evalii_{n-\dimevec}, \\
\delta_{2,+} &=  \evalii_{n-\dimevec+1} - c_1 \evali_{n-\dimevec} - c_0,  \\
\delta_{2,-} &=  c_1 \evali_{n-\dimevec} + c_0 - \evalii_{n}.
\end{align*}

The above arguments  can be applied to these values of $\delta$ without complications. Consequently, for $j=n-r,$ the cost function in \eqnref{eqn_affine_general_bound}  tends to $1$ as $c_0$ tends to either $\infty$ or $-\infty.$

\end{document}